\documentclass{svproc}
\usepackage{url}

\usepackage{amssymb}
\usepackage{amsmath}
\usepackage{comment}
\usepackage{tikz}

\tikzset{vtx/.style={inner sep=1.7pt, outer sep=0pt, circle, fill,draw}}

\tikzset{vtxB/.style={inner sep=1.7pt, outer sep=0pt, circle, fill=white,draw}}

\author{J\'ozsef Balogh\inst{1} \and Felix Christian Clemen\inst{2}
 \and Bernard Lidick\'{y}\inst{3}}

\institute{Department of Mathematics, University of Illinois at Urbana-Champaign, Urbana, Illinois 61801, USA,
\email{jobal@illinois.edu}
\and
Department of Mathematics, University of Illinois at Urbana-Champaign, Urbana, Illinois 61801, USA,
\email{fclemen2@illinois.edu}
\and 
Iowa State University, Department of Mathematics, Iowa State University, Ames, Iowa, USA 
\email{lidicky@iastate.edu}}

\title{Max Cuts in Triangle-free Graphs }
\titlerunning{Erd\H{o}s' Conjecture on Cuts in Triangle-free Graphs}

\begin{document}
\mainmatter

\maketitle 
\begin{abstract}
A well-known conjecture by Erd\H{o}s states that every triangle-free graph on $n$ vertices can be made bipartite by removing at most $n^2/25$ edges. This conjecture was known for graphs with edge density at least $0.4$ and edge density at most $0.172$. Here, we will extend the edge density for which this conjecture is true; we prove the conjecture for graphs with edge density at most $0.2486$ and for graphs with edge density at least $0.3197$. Further, we prove that every triangle-free graph can be made bipartite by removing at most $n^2/23.5$ edges improving the previously best bound of $n^2/18$.   
\end{abstract}

\section{Introduction}
How many edges need to be removed from a triangle-free graph on $n$ vertices to make it bipartite? Erd\H{o}s~\cite{MR0409246} asked this question and conjectured that $n^2/25$ edges would always be sufficient. This would be sharp as the balanced blow-up of $C_5$ with class sizes $n/5$ needs at least $n^2/25$ edges removed to be made bipartite. For a graph $G$, denote $D_2(G)$ the minimum number of edges which have to be removed to make $G$ bipartite.
\begin{conjecture}{(Erd\H{o}s\cite{MR0409246})}
\label{conj:main}
For every triangle-free graph $G$ on $n$ vertices 
\begin{equation}
    D_2(G)\leq \frac{n^2}{25}.
\end{equation}
\end{conjecture}
An elementary probabilistic argument (see e.g.~\cite{Howmany}) resolves Conjecture~\ref{conj:main} for graphs $G$ with at most $2/25 n^2$ edges: Take a random bipartition where each vertex, independently from each other, is placed  with probability $1/2$ in one of the two classes. The expected number of edges inside both of the classes is $|E(G)|/2$. Thus, there exists a bipartition with at most $|E(G)|/2$ edges inside the classes. Note that this argument does not use that $G$ is triangle-free. Erd\H{o}s, Faudree, Pach and Spencer~\cite{Howtotriangle} slightly improved this random cut argument utilizing triangle-freeness.   
\begin{theorem}[Erd\H{o}s, Faudree, Pach, Spencer~\cite{Howtotriangle}]
\label{ErFaPaSp}
For every triangle-free graph with $n$ vertices and $m$ edges
\begin{equation}
    D_2(G)\leq \min \left\{\frac{m}{2}-\frac{2m(2m^2-n^3)}{n^2(n^2-2m)},m-\frac{4m^2}{n^2}   \right\} \leq \frac{n^2}{18}.
\end{equation}
\end{theorem}
This confirmed Conjecture~\ref{conj:main} for graphs with roughly at most $0.086 n^2$ edges and graphs with at least $n^2/5$ edges. It also gives the current best bound on the Erd\H{o}s problem; one can remove at most $n^2/18$ edges to make a triangle-free graph bipartite. We improve this result and extend the range for which Erd\H{o}s' conjecture is true. 
\begin{theorem}
\label{improvedtriangle}
Let $G$ be a triangle-free graph on $n$ vertices. Then, for $n$ large enough,
\begin{itemize}
    \item[(a)] $D_2(G)\leq \frac{n^2}{23.5},$
    \vspace{0.1cm}
    \item[(b)] $D_2(G) \leq \frac{n^2}{25}$ when  $|E(G)|\geq 0.3197 \binom{n}{2}$,  % 0.31960253 
    \vspace{0.1cm}
   \item[(c)] $D_2(G) \leq \frac{n^2}{25}$ when   $|E(G)|\leq 0.2486 % 0.24869221
    \binom{n}{2}$.
\end{itemize}
\end{theorem}
Sudakov studied a related question; he~\cite{MR2359832} determined the maximum number $D_2(G)$ for $K_4$-free graph $G$. Recently, Hu, Lidick\'{y}, Martins, Norin and Volec~\cite{hulimanovo} announced a proof for determining the maximum number $D_2(G)$ for $n$-vertex $K_6$-free graphs $G$.
They use the method of flag algebras, developed by Razborov~\cite{flagsRaz}, to describe local cuts which leads to the solution. We use a similar idea of encoding local cuts.

Our proof of Theorem~\ref{improvedtriangle} also extends on the ideas from Erd\H{o}s, Faudree, Pach, Spencer~\cite{Howtotriangle}. While their proof uses two different ways of finding bipartitions, our proof uses many ways. In order to handle a large amount of bipartitions, we use the method of flag algebras. It relies on formulating a problem as a semidefinite program and then using a computer to solve it.

We will handle graphs with edge density close to $2/5$ (the density of the conjectured extremal example) separately. In this range we use standard techniques from extremal combinatorics, such as a minimum degree removing algorithm. Additionally, we will make use of the following result by Erd\H{o}s, Gy\H{o}ri and Simonovits~\cite{Howmany}.

A $C_5$-blow-up $H$ is a graph with vertex set $V(H)=A_1 \cup A_2 \cup A_3 \cup A_4 \cup A_5$ and edges $xy\in E(H)$ iff $x\in A_i$ and $y\in A_{i+1}$ for some $i\in [5]$, where $A_6:=A_1$. 
\begin{theorem} [Erd\H{o}s, Gy\H{o}ri and Simonovits~\cite{Howmany}]
\label{ErdosC_5blowup}
Let $G$ be a $K_3$-free graph on $n$ vertices with at least $n^2/5$ edges. Then there exists an unbalanced blow-up of $C_5$ $H$ such that   
\begin{equation}
D_2(G)\leq D_2(H).
\end{equation} 
\end{theorem}
Note that this result recently was extended to cliques by Kor\'andi, Roberts and Scott~\cite{scott} confirming a conjecture from~\cite{Makingrpartite}.

There is a local version of Conjecture~\ref{conj:main}.
\begin{conjecture}{(Erd\H{o}s\cite{MR0409246})}
\label{conj:second}
Every triangle-free graph on $n$ vertices contains a vertex set of size $\lfloor n/2 \rfloor$ that spans at most $n^2/50$ edges.
\end{conjecture}
Erd\H{o}s~\cite{MR1439273} offered \$250 for the first solution of this conjecture. As pointed out by Krivelevich~\cite{MR1320169}, for regular graphs Conjecture~\ref{conj:second} would imply Conjecture~\ref{conj:main}. We are wondering if similar methods we are using could be used to make progress towards proving Conjecture~\ref{conj:second}.

This extended abstract is organized as follows. In Section~\ref{flagsection} we present our setup for flag algebras to give a sketch of the proof of the main part of Theorem~\ref{improvedtriangle}. In Section~\ref{slighlty below} we sketch the proof of Conjecture~\ref{conj:main} in the edge range slightly below edge density $2/5$. 
%Finally, in Section~\ref{sec:furtherdiscussion} we %discuss how one could proceed making progress on this %problem.  

\section{Proof Sketch of Theorem~\ref{improvedtriangle}}
\subsection{Setup for flag algebras}
\label{flagsection}
\begin{comment}
We start by introducing some notation. For a graph $G$ with vertex set $V(G)$ and edge set $E(G)$, we denote with $e(A)$ the number of edges in the induced graph $G[A]$. For disjoint sets $A,B\subset V(G)$ we denote $e(A,B)$ the number of crossing edges between $A$ and $B$, i.e. all edges with one endpoint in $A$ and the endpoint in $B$. Further, $e^c(A,B)$ denotes the number of non edges between $A$ and $B$, that is $e^c(A,B)=|A||B|-e(A,B)$.

Given a graph $G$ with vertices $v_1,...,v_n$ and integers $j_1,...,j_n$, the graph $G[j_1,j_2,\ldots,j_n]$ (\textupt{blow-up} of $G$) is defined as follows: The $i$-th vertex of $G$ $v_i$ is replaced by $j_i$ new vertices $v_{i,1},\ldots v_{i,j_i}$. A vertex $v_{i,j}$ is adjacent to $v_{i',j'}$ in $G[j_1,j_2,\ldots,j_n]$ iff $v_i$ is adjacent to $v_{i'}$ in $G$. 
\end{comment}

Towards contradiction assume that there is a triangle-free graph $G$ on $n$ vertices with $D_2(G)\geq n^2/25$.
This means that whenever we create a bipartition of $V(G)$, then it has at least $n^2/25$ edges inside the two parts.
Using flag algebras, one can define bipartitions and count edges inside of the two parts.

\newcommand{\vc}[1]{\ensuremath{\vcenter{\hbox{#1}}}}
\tikzset{vtx/.style={inner sep=1.7pt, outer sep=0pt, circle, fill}} 
\tikzset{unlabeled_vertex/.style={inner sep=1.7pt, outer sep=0pt, circle, fill}} 
\tikzset{labeled_vertex/.style={inner sep=2.2pt, outer sep=0pt, rectangle, fill=yellow, draw=black}} 
\tikzset{edge_color0/.style={color=black,line width=1.2pt}} 
\tikzset{edge_color1/.style={color=red,  line width=1.2pt,opacity=0}} 
\tikzset{edge_color2/.style={color=blue, line width=1.2pt,opacity=1}}

For example, in a graph $G$ one could fix a vertex $v$ and define the bipartition of $G$ as $V(G) = N(v) \cup \overline{N(v)}$.
If one uses this bipartition, all edges in $\overline{N(v)}$ need to be removed while $N(v)$ is independent since $G$ is triangle-free.
This can be written in flag algebras in the following way
\begin{align}\label{eq:v}
\vc{\begin{tikzpicture}[scale=0.5]
\draw \foreach \x in {0,1,2}{(270+\x*360/3:0.8) coordinate(x\x)};
\draw[edge_color2] (x1)--(x2);
\draw (x0) node[labeled_vertex,label=below:$v$]{};
\draw (x1) node[unlabeled_vertex]{};
\draw (x2) node[unlabeled_vertex]{};
\end{tikzpicture}}
\geq \frac{2}{25},
\end{align}
where the depicted graph represents its expected induced density when unordered pair of black vertices is picked uniformly at random while the yellow vertex is fixed.
In proving Theorem~\ref{ErFaPaSp}, Erd\H{o}s, Faudree, Pach, Spencer~\cite{Howtotriangle} used this cut and the following cut.
Let $uv$ be two adjacent vertices. Let $N(u)$ be one part and $N(v)$ be the other part. The remaining vertices in $\overline{N(u)\cup N(v)}$ are partitioned uniformly at random with probability $1/2$ to either of the two parts.
Since $G$ is $K_3$-free, one obtains the following equation for flag algebras
\begin{align}\label{eq:K2}
\frac{1}{2}
\,\,
\vc{\begin{tikzpicture}[scale=0.5]
\draw \foreach \x in {0,1,2,3}{(270-45+\x*360/4:0.8) coordinate(x\x)};
\draw[edge_color2] (x0)--(x1);
\draw[edge_color2] (x3)--(x2);
\draw[edge_color2] (x0)--(x3);
\draw (x0) node[labeled_vertex,label=below:$u$]{};
\draw (x1) node[labeled_vertex,label=below:$v$]{};
\draw (x2) node[unlabeled_vertex]{};
\draw (x3) node[unlabeled_vertex]{};
\end{tikzpicture}} 
+
\frac{1}{2}
\,\,
\vc{\begin{tikzpicture}[scale=0.5]
\draw \foreach \x in {0,1,2,3}{(270-45+\x*360/4:0.8) coordinate(x\x)};
\draw[edge_color2] (x0)--(x1);
\draw[edge_color2] (x3)--(x2);
\draw (x0) node[labeled_vertex,label=below:$u$]{};
\draw (x1) node[labeled_vertex,label=below:$v$]{};
\draw (x2) node[unlabeled_vertex]{};
\draw (x3) node[unlabeled_vertex]{};
\end{tikzpicture}} 
+
\frac{1}{2}
\,\,
\vc{\begin{tikzpicture}[scale=0.5]
\draw \foreach \x in {0,1,2,3}{(270-45+\x*360/4:0.8) coordinate(x\x)};
\draw[edge_color2] (x0)--(x1);
\draw[edge_color2] (x1)--(x2);
\draw[edge_color2] (x3)--(x2);
\draw (x0) node[labeled_vertex,label=below:$u$]{};
\draw (x1) node[labeled_vertex,label=below:$v$]{};
\draw (x2) node[unlabeled_vertex]{};
\draw (x3) node[unlabeled_vertex]{};
\end{tikzpicture}} 
\geq \frac{2}{25}.
\end{align}

This idea of defining cuts can be generalized by rooting on more vertices. 
Pick a copy of a \emph{labeled} graph $H$ on $k$ vertices in $G$. 
This will partition the rest of $V(G)$ into classes $X_1,\ldots,X_{2^k}$ based on the adjacencies to the fixed $k$ vertices.
Now we construct a bipartition of $V(G)$ into sets $A$ and $B$. For each class $X_i$ fix $p_i \in [0,1]$ and for each vertex in $X_i$ 
we put it to $A$ with probability $p_i$ and to put it to $B$ otherwise, i.e., with probability $(1-p_i)$.

This creates a bipartition and it is possible to count the edges that need to be removed using flag algebras. We can include all cuts rooted on at most $4$ vertices and $C_5$.
\begin{enumerate}
\item $|V(H)| \leq 2$ and $p_i \in \{0,0.5,1\}$, gives 10 cuts,
\item $|V(H)| \leq 3$ and $p_i \in \{0,0.5,1\}$, gives 108 cuts,
\item $|V(H)| = 4$ and $p_i \in \{0,1\}$, gives 953 cuts,
\item $H = C_5$,  and $p_i \in \{0,1\}$, gives 125 cuts.
%\item $|V(H)| = 5$, $|V(E)| \geq 3$, and $p_i \in \{0,1\}$, gives 8313 cuts.
\end{enumerate}
However, for $k\geq 6$, there are more possible inequalities than computers can reasonably handle. Therefore we have to decide on which we want to use. We will present two particular important ones here. 
\begin{comment}
One could then add constraints that for any cut that is used, one has to remove at least $1/25 n^2$ edges and see what else happens.
First thing we need to do is to establish what types of cuts are being used.
We have all cuts, where
\end{comment}
%Jozis note: about 0.5 is useful only for the block that is not incident to anything - if you split, it is a linear function, only thing that matters is that if it is a neighborhood, it is independent set - so we get nothing about adding it.
%Adding all cuts actually slows down the computation quite a bit. Maybe the program is not efficient.

%Using these cuts, we obtain $D_2(G) \leq n^2/25$ if $|E(G)| \geq 0.351n$.

%Say upper bound, cuts 1.--5. give 0.3504.
%\end{comment}

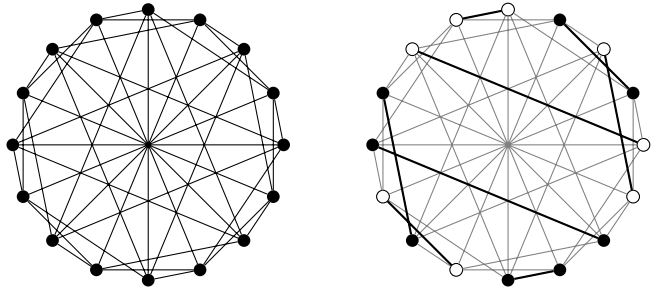
\begin{figure}
\begin{center}
    \begin{tikzpicture}[scale=0.6]
\draw 
\foreach \x in {0,...,15}{(22.5*\x:3)node[vtx](\x){}}
\foreach \x/\y in {0/1,0/6,0/8,0/10,0/13,1/3,1/4,1/9,1/15,2/3,2/8,2/10,2/12,2/15,3/5,3/6,3/11,4/5,4/10,4/12,4/14,5/7,5/8,5/13,6/7,6/12,6/14,7/9,7/10,7/15,8/9,8/14,9/11,9/12,10/11,11/13,11/14,12/13,13/15,14/15}{
(\x)--(\y)
}
;
    \end{tikzpicture}
    \hspace{0.8cm}
    \begin{tikzpicture}[scale=0.6]
\draw 
\foreach \x in {0,...,15}{(22.5*\x:3)node[vtx](\x){}}
\foreach \x in {0,2,4,5,6,9,11,15}{(22.5*\x:3)node[vtxB](\x){}}
;
\draw[color=gray]
\foreach \x/\y in {0/1,0/6,0/8,0/10,0/13,1/3,1/4,1/9,1/15,2/3,2/8,2/10,2/12,2/15,3/5,3/6,3/11,4/5,4/10,4/12,4/14,5/7,5/8,5/13,6/7,6/12,6/14,7/9,7/10,7/15,8/9,8/14,9/11,9/12,10/11,11/13,11/14,12/13,13/15,14/15}{
(\x)--(\y)
}
;
\draw[thick]
\foreach \x/\y in {0/6,1/3,2/15,4/5,7/10,8/14,9/11,12/13}{
(\x)--(\y)
}
;

;
    \end{tikzpicture}    
\end{center}
\caption{Clebsch graph and its cutting}\label{fig:clebsch}
\end{figure}

Norin and Ru Sun~\cite{norin2016triangleindependent} observed that the Clebsch graph, see Figure~\ref{fig:clebsch}, is particularly unfriendly when applying local cuts. We add cuts that are specially designed to cut the Clebsch graph. %, which yelds the improvement to 0.3472.
The root is a 4-cycle $v_0v_1v_2v_3v_0$ and two additional vertices $v_4$ and $v_5$ with edges $v_4v_0$ and $v_1v_5$. Although this is a bipartite graph, we create a bipartition as if $v_1,v_2,v_5$ and $v_1,v_3,v_4$ were in the same parts respectively.                                        

\begin{comment}
It has 40 edges and $D_2(CG)=8$. Conjecture~\ref{conj:main} asks for $D_2(CG) \leq 16^2/25 = 10.24$.
We use a cut 
\begin{verbatim}
    4   5
    |   |
    0---1
    |   |
    |   |
    2---3

    0,2,5 are one part    -> all neighbors only 1 or only 3 + N(0,3) & not 5 + N(1,2,4)
    1,3,4 is another part -> all neighbors only 0, 2 + N(0,5,3) + N(1,2) & not 4
\end{verbatim}
%    8 6  2 2 1 2 1 0 0  1 2 1 2 0 0  2 1 1 0 0  1 1 0 0  1 0 0   0 0  2
\end{comment}    

Another inequality that made a big difference is an extension of \eqref{eq:K2}.
While \eqref{eq:K2} partitions neighbors of the chosen two vertices very well,
the non-neighbors can be partitioned better.
In particular, we pick another $K_2$ in the non-neighborhood and do the same partition once more. This results in rooting on $K_2 \cup K_2 \cup K_2$.

\begin{comment}
\begin{theorem}
\label{flagErdos}
There exists $n_0$ such that for all $n\geq n_0$ the following holds. Let $G$ be an $n$-vertex graph with $|E(G)|\leq (0.2-10^{-8})n^2$ edges, then
$D_2(G)\leq n^2/25$.
\end{theorem}
\end{comment}
Our flag algebra proof cannot deal with the density range close to $2/5$, i.e. close to the conjectured extremal example. In the following section we explain how this density range can be handled.  

\subsection{High density range}
\label{slighlty below}

In this section we provide a sketch of the proof of Erd\H{o}s' conjecture for graphs with edge density slightly below $2/5$. 
\begin{theorem}
\label{theoslightlybelow}
There exists $n_0$ such that for all $n\geq n_0$ the following holds. Let $G$ be an $n$-vertex triangle-free graph with $ |E(G)|\geq (0.2-\varepsilon) n^2$ edges, where $\varepsilon=10^{-8}$. Then $D_2(G)\leq n^2/25$.  
\end{theorem}
Let $G_n:=G$ be a triangle-free graph on $n$ vertices with $|E(G)|\geq (0.2-\varepsilon) n^2$ edges. Assume, towards contradiction, $D_2(G)>n^2/25$. We iteratively remove a vertex of minimum degree from $G$. This means $G_i=G_{i+1}-x$, where $\deg(x)=\delta(G_{i+1})$. We stop this algorithm if $\delta(G_i)> \frac{3}{8}i$ or after $\lfloor 5 \cdot 10^{-7}n\rfloor$ rounds. Let $m$ be the stage in which the algorithm stops. 
\begin{lemma}
\label{cutclaim}
We have
\begin{align}
    D_2(G)\leq \frac{3}{32}(n^2-m^2+n-m) +D_2(G_m).
\end{align}
\end{lemma}
This Lemma can be verified by taking a smallest cut of $G_m$ and adding the remaining vertices to the set where they have smaller neighborhood in. 

Depending on when the algorithm stops we perform a different analysis. If the algorithm stops ``late'', then $G_m$ has edge density of slightly more than $2/5$.
By Lemma~\ref{cutclaim}, we can assume that 
\begin{equation}
    D_2(G_m) \geq \frac{n^2}{25}-\frac{3}{32}(n^2-m^2+n-m).
\end{equation}
By Theorem~\ref{ErdosC_5blowup} we can find a $C_5$-blow-up $H$ on $m$ vertices with classes $A_1,A_2,A_3,$ $A_4,A_5$ satisfying $|E(H)|\geq |E(G_m)|$ and $D_2(H)\geq D_2(G_m)$. In fact, it can also be assumed that the class sizes of $H$ are symmetric, that is $|A_2|=|A_5|+o(n)$ and $|A_3|=|A_4|+o(n)$. A straight-forward optimization of the number of edges in $H$ gives a contradiction with $|E(H)|\geq |E(G_m)|$.

If the algorithm stops early, we make use of a result by H{\"a}ggkvist~\cite{MR671908} who proved that every triangle-free graph on $n$ vertices with minimum degree more than $3n/8$ is a subgraph of a $C_5$-blow-up. Having this particular structure, it can be calculated that $D_2(G)\leq n^2/25$, we omit the detailed computations.

\subsection{Concluding Remarks}
Note that Theorem~\ref{improvedtriangle} only holds for $n\geq n_0$ for some $n_0$ large enough. However, this is not an actual restriction towards proving Conjecture~\ref{conj:main}. Assuming Conjecture~\ref{conj:main} were to hold for all $n\geq n_0$, then it actually holds for all $n$ by the following argument. Let $G$ be a triangle-free graph on $n <n_0$ vertices and assume, towards contradiction, that $D_2(G)>n^2/25$. Consider the blow-up $G'$ of $G$, where each vertex is replaced by an independent set of size $\lceil \frac{n_0}{n}\rceil$ and two vertices in different sets are made adjacent iff the corresponding vertices in $G$ were adjacent. This new graph $G'$ is still triangle-free and has at least $n_0$ vertices. A result by Erd\H{o}s, Gy\H{o}ri and Simonovits~\cite[Theorem 7]{Howmany} gives 
\begin{equation}
    \frac{D_2(G')}{\left(\lceil \frac{n_0}{n}\rceil n\right)^2}\geq \frac{D_2(G)}{n^2} >\frac{1}{25},
\end{equation}
contradicting that we assumed Conjecture~\ref{conj:main} holds for all $n\geq n_0$ and therefore in particular for $G'$.

We believe Theorem~\ref{improvedtriangle} can be improved by adding more cuts to the calculation and possibly lead to the proof of Conjecture~\ref{conj:main}. Adding more cuts lead to marginal improvements so far. We are looking at other cuts as well but the time needed to perform the calculations grows quickly and it may take a while until a significant improvement is obtained. 

\subsection{Acknowledgements}
We thank Humberto Naves, Florian Pfender, and Jan Volec for fruitful discussions in early stages of this project. 
The first author is supported by NSF grants DMS-1764123, DMS-1937241, the Langan Scholar Fund (UIUC), and the Simons Fellowship. The first and second authors are supported by the Arnold O. Beckman Research
Award (UIUC RB 18132).
The last author is supported by NSF grant DMS-1855653.

%This is an extended abstract for Eurocomb 2021. Comments are welcome.

\bibliographystyle{abbrv}
\bibliography{bibilo}

\end{document}